\documentclass[10pt,twoside]{article}
\usepackage{graphicx}
\usepackage{amsmath}
\usepackage{Latex-document}
\usepackage{amssymb}
\title{\bf  Equivariant Bloch-Kato conjecture and
non-abelian Iwasawa Main Conjecture}

\author{A. Huber\thanks{Math. Institut, Universit{\"a}t 
Leipzig, Augustusplatz 10/11, 04109 Leipzig, Germany.
E-mail: huber@mathematik.uni-leipzig.de} {\ and} 
G. Kings\thanks{
NWF I-Mathematik, Universit{\"a}t Regensburg, 93040 Regensburg Germany.
E-mail: guido.kings@mathematik.uni-regensburg.de
}}

\date{June 2002}

\newcommand{\Q}{{\Bbb{Q}}}  
\newcommand{\C}{{\Bbb{C}}}  
\newcommand{\Z}{{\Bbb{Z}}}  
\newcommand{\R}{{\Bbb{R}}} 
\newcommand{\1}{\mathbf{1}}
\newcommand{\Dh}{\mathcal{D}}
\newcommand{\Oh}{\mathcal{O}}
\newcommand{\Th}{\mathcal{T}}
\newcommand{\Lh}{\mathcal{L}}

\newcommand{\tensor}{\otimes} 
\newcommand{\isom}{\cong} 
\newcommand{\ord}{\operatorname{ord}} 
\newcommand{\Aut}{\operatorname{Aut}}
\newcommand{\Gl}{\mathrm{Gl}}
\newcommand{\bild}{\operatorname{Im}}
\newcommand{\Hom}{\operatorname{Hom}}
\newcommand{\End}{\operatorname{End}}

\renewcommand{\det}{\operatorname{det}}
\newcommand{\rn}{\mathrm{rn}}
\newcommand{\deltatilde}{\widetilde{\delta}}

\newcommand{\ab}{\mathrm{ab}}
\newcommand{\triv}{\mathrm{triv}}
\newcommand{\hm}{H_{\mathcal{M}}}
\newcommand{\prolim}{\varprojlim}
\newcommand{\Gal}{\operatorname{Gal}}

 \newcommand{\rem}{\noindent{\bf Remark:\ }}     
\newtheorem{conj}{Conjecture}[subsection]
\newtheorem{prop}[conj]{Proposition}
\newtheorem{cor}[conj]{Corollary}
\newtheorem{lemma}[conj]{Lemma}
\newtheorem{thm}[conj]{Theorem}
\begin{document}
\maketitle

\thispagestyle{first} \setcounter{page}{1}

\begin{abstract}
In this talk we explain the relation between the (equivariant) Bloch-Kato
conjecture for special values of $L$-functions and the Main Conjecture 
of (non-abelian) Iwasawa theory. On the way we will discuss briefly the case of Dirichlet characters in
the abelian case. We will also discuss how ``twisting'' in the non-abelian
case would allow to reduce the general conjecture to the case of number fields.
This is one the main motivations for a non-abelian Main Conjecture.


\end{abstract}

\vskip 12mm

\section{Introduction} 
The class number formula expresses the leading coefficient of a Dedekind-$\zeta$-function of a number field $F$ in terms of arithmetic invariants of $F$:
\[ \zeta_F(0)^*= -\frac{h R_F}{w_F} \]
($h$ the class number, $R_F$ the regulator, $w_F$ the number of roots of
unity in $F$). By work of Lichtenbaum, Bloch, Beilinson, and Kato  among others,
the class number formula has
been generalized to other $L$-functions of varieties
(or even motives) culminating in the
Tamagawa number conjecture by Bloch and Kato.

Iwasawa, on the other hand, initiated the study of the growth of the class numbers in towers
of number fields. His decisive idea was to
consider the class group of the tower as a module under the completed
group ring of the Galois group of the tower.
From his work evolved the ``Main Conjecture'' describing
this growth in terms of the $p$-adic $L$-function. 
 
It is a surprising insight of Kato that an equivariant version of the Tamagawa number
conjecture can be viewed as a version of the Main Conjecture of Iwasawa theory.
Perrin-Riou, in her efforts to develop a theory of $p$-adic $L$-functions, arrived
at a similar conclusion.

The purpose of this paper is to make the connection between the equivariant Tamagawa number
conjecture and the Iwasawa Main Conjecture precise. In the spirit of Kato, we formulate 
an Iwasawa Main Conjecture (\ref{mainconj}) for arbitrary motives and towers of number fields whose
Galois group is a $p$-adic Lie group. This formulation does not involve $p$-adic $L$-functions. We
show that it is implied by the equivariant Tamagawa number conjecture formulated by Burns and Flach.
For ease of exposition, we restrict to the case of $L$-values at very negative integers, where the
Bloch-Kato exponential does not play a role. 
The study of non-abelian Iwasawa theory was initiated by Coates. Recently,
there have been systematic results by Coates, Howson, Ochi, Schneider, Sujatha and Venjakob.

Our interest in allowing general towers  of number fields is motivated by the possibility 
of reducing the Tamagawa number conjecture to an equivariant class number formula (modulo hard
conjectures, see \ref{strategy}).

Important special cases of the Main Conjecture were considered by (alphabetical order) Coates, Greenberg,
Iwasawa, Kato, Mazur, Perrin-Riou, Rubin, Schneider, Wiles and more recently by Ritter and Weiss.

It is a pleasure to thank C. Deninger, S. Howson, B. Perrin-Riou, A. Schmidt, P. Schneider for
helpful comments and discussions.

\section{Non-abelian equivariant Tamagawa number conjecture}
\subsection{Notation}\label{notation}
Fix $p\neq 2$ and let $M$ be a motive over $\Q$ with coefficients in $\Q$, for example $M=h^r(X)$, $X$ a smooth
projective variety over $\Q$. It has Betti-realization $M_B$ and $p$-adic
realization $M_p$.
Let $M^\lor$ be the dual motive. In the $p$-adic realization it corresponds to the dual Galois module.
We denote by $\hm^1(\Z, M(k))$ the ``integral'' motivic cohomology of the motive $M$ in the
sense of Beilinson \cite{Beilinson}. 

For any finite Galois extension $K/\Q$ with Galois group $G$, let $\Q[G]$ be
the group ring of $G$. It is a non-commutative ring with center denoted $Z(\Q[G])$.

We consider the {\em deformation}
$\Q[G]\tensor M:= h^0(K)\tensor M$. 
If $M=h^r(X)$ and $K/ \Q$ is a number field, then $h^0(K)\tensor M=h^r(X\times K)$ considered
as a motive over $\Q$.

We consider a finite set of primes $S$ satisfying:
\begin{description}
\item[(*)]
$\Q[G]\otimes M$ and $K$ have good  reduction at all primes not dividing $S$, and 
$p \in S$.
\end{description}
\subsection{Equivariant $L$-functions}
We {\em assume} the usual conjectures about the $L$-functions of
motives, like meromorphic continuation and functional equation etc., to be satisfied.

In order to define the {\em equivariant}
$L$-function for $G$ (without the Euler factors at the primes dividing $S$), 
consider a Galois extension $E/ \Q$ such that $E[G]\isom \bigoplus_\rho \End_E(V(\rho))$, where
$V(\rho)$ are absolutely irreducible representations of $G$.
Then the center of $E[G]$ is $Z(E[G])\isom \bigoplus_\rho E$ and the  motives $V(\rho)\tensor M$ 
have coefficients in  $E$. 
We  define 
\[
L_S(G, M, k)^*:= \Big( L_S(V(\rho)\tensor M, k)^*\Big)_\rho \in Z(E\otimes_\Q\C[G])^*
\]
to be the element with $\rho$-component the leading coefficient at $s=k$ of the $E\otimes_\Q \C$-valued  
$L$-functions $L_S(V(\rho)\otimes M, s)$ without the 
Euler factors at $S$. Then $L_S(G, M, k)^*$ has actually values in $Z(\R[G])^*$ 
(see \cite{BF} Lemma 7) and is 
independent of the choice of $E$. We will always consider 
$L_S(G, M, k)^*$ as an element in $Z(\R[G])\subset \R[G]$. 

\rem In  \cite{Kato-mc} Kato uses a different description of this equivariant $L$-function.

\subsection{Non-commutative determinants}
We follow the point of view of Burns and Flach.
Let $A$ be a (possibly non-commutative) ring and $V(A)$ the category
of virtual objects in the sense of Deligne \cite{Deligne}. $V(A)$ is a monoidal tensor category
and
has a unit object $\1_A$. Moreover it is a groupoid, i.e., all
morphisms are isomorphisms. There  is a functor
\[
\det_A: \{ \text{perfect complexes of $A$-modules and isomorphisms}\} \to V(A)
\]
which is 
multiplicative on short
exact sequences. 
The group of isomorphism classes
of objects of $V(A)$ is $K_0(A)$ and
\[ \Aut(\1_A)=K_1(A)=\Gl_\infty(A)/E(A) \]
($E(A)$ the elementary matrices). In general $\Hom_{V(A)}(\det_A X,\det_A Y)$
is either empty or a $K_1(A)$-torsor. 
If $A\to B$ is a ring homomorphism, we get a functor 
$B \tensor :V(A)\to V(B)$ such that tensor product and $\det_A$ commute.

\noindent{\bf Convention:}
By abuse of notation we are going
to write $z\in \det_A X$ for $z:\1_A\to \det_A X$ and call this a
{\em generator} of $\det_A X$.

If $A$ is commutative and local, then the category of virtual objects is
equivalent to the category of pairs $(L,r)$ where $L$ is an invertible
$A$-module and $r\in \Z$.
One recovers the theory of determinants of Knudson and Mumford.
The unit object is $\1_A=(A,0)$ and one has $ \Aut(\1_A)=K_1(A)=A^*$.
Thus $K_1(A)$ is used as  generalization of 
$A^*$ to the non-commutative case. Generators of $\det_A X=(L,0)$ in the above
sense correspond to $A$-generators of $L$.

\subsection{Formulation of the conjecture}\label{2.4}
The original conjecture dates back to Beilinson \cite{Beilinson} and 
Bloch-Kato \cite{Bloch-Kato}. The idea of an equivariant formulation
is due to Kato \cite{KatoGalaxy} and \cite{Kato-mc}. Fontaine and
Perrin-Riou gave a uniform formulation for mixed motives and all
values of $L$-functions at all integer values
\cite{FontaineBourbak}, \cite{Fo-PR}. The generalization to non-abelian
coefficients is due to Burns and Flach \cite{BF}.

For simplicity of
exposition, we restrict to values at very negative integers. In the
absolute case this coincides with the formulation given by Kato
in \cite{KatoGalaxy}. We consider a motive $M$ and values at $1-k$ where
$k$ is {\em big enough}. 
In the case $M=h^r(X)$, $k$ big enough means that 
\begin{itemize}
\item
$k>\inf \{r, \dim (X)\} $, $(r,k)\neq (1,0);(2\dim(X),\dim(X)+1)$ and $2k\neq r+1$.
\item for all $ \ell \in S$ the local Euler factor  $L_\ell(M_p^\lor,s)^{-1}$  at $\ell $ 
does not vanish at $1-k$.
\end{itemize}

Consider the (injective) reduced norm map $rn: K_1(\R[G])\to Z(\R[G])^*$ and
recall that $L_S(G, M^\lor, 1-k)^*\in Z(\R[G])^*$. By strong approximation
(see \cite{BF} Lemma 8) there is $\lambda\in Z(\Q[G])^*$ such
that $\lambda L_S(G,M^\lor (1-k))^*$ is in the image of $K_1(\R[G])$ under $rn$.
Let
\[ \lambda L_S(G,M^\lor (1-k))^* \in \1_{\R[G_n]} \]
be the corresponding generator.
For $k$ big enough, we define the {\em fundamental line} in $V(\Q[G])$ as
\[ \Delta_f(G,M^\lor(1-k))=\det^{-1}_{\Q[G]} \hm^1(\Z,\Q[G]\tensor M(k))
  \tensor \det_{\Q[G]} (\Q[G]\tensor M_B(k-1))^+ \ .
\]
Here $+$ denotes the fixed part under complex conjugation.

\begin{conj}\label{conj}
Let $M$ be as in \ref{notation}, $p\neq 2$ a prime and $k$ big enough.
\begin{enumerate}
\item The Beilinson regulator $r_{\Dh}$ induces an isomorphism
\[ \Delta_f(G,M^\lor(1-k))\tensor \R\isom \1_{\R[G]} \]
\item Under this isomorphism the generator $(\lambda L_S(G,M^\lor (1-k))^*)^{-1}$ is 
induced by a
(unique) generator
\[ (\lambda^{-1} \delta(G,M,k)) \in\Delta_f(G,M^\lor(1-k)) \]
\end{enumerate}
The reduced norm is an isomorphism $K_1(\Q_p[G])\isom Z(\Q_p[G])^*$. 
Using the operation of $K_1(\Q_p[G])$ on generators in
$ \Delta_f(G,M^\lor(1-k))\tensor\Q_p$, we put
\[ \delta_p(G,M,k):=(\lambda^{-1} \delta(G,M,k) )\lambda\in \Delta_f(G,M^\lor(1-k))\tensor\Q_p \]
Note that this generator is independent of the choice of $\lambda$.
\begin{enumerate}
\setcounter{enumi}{2}
\item The $p$-adic regulator $r_p$ induces an isomorphism
\begin{multline*}
\Delta_f(G,M^\lor(1-k))\tensor \Q_p\isom\\ 
   \det_{\Q_p[G]}^{-1} H^1(\Z[1/S], \Q_p[G]\tensor M_p(k)) \tensor
   \det_{\Q_p[G]} (\Q_p[G]\tensor M_B(k-1))^+
\end{multline*}
\item Let $T_B\subset M_B$ be a lattice such that
$T_p=T_B\tensor \Z_p\subset M_p$ is Galois stable.
Under the last isomorphism  $\delta_p(G,M,k)$ is induced by a generator
\[ \deltatilde_p(G,M,k)\in \det_{\Z_p[G]} R\Gamma(\Z[1/S],\Z_p[G]\tensor T_p(k)) \tensor
   \det_{\Z_p[G]} ( \Z_p[G]\tensor T_B)(k-1))^+
\]
\end{enumerate}
\end{conj}

\rem
a) The conjecture is compatible with change of group $G$. If $G\to G'$ is
a surjection, then the equivariant conjecture for $G$ tensored with
$\Q[G']$ over $\Q[G]$ gives the conjecture for $G'$.\\
b) The element $\deltatilde_p(G,M,k)$ is determined up to an element in
the kernel of the map $K_1(\Z_p[G])\to K_1(\Q_p[G])$. In the commutative case,
this map is always injective.\\
c) The conjecture is independent of $T$. It is also independent of $S$. This
computation shows that the definition of the equivariant $L$-function forces the use
of the reduced norm in the formulation of the conjecture.

\section{Non-abelian Main Conjecture}
\subsection{Iwasawa algebra and modules}\label{3-1}
Let $K_n$ be a tower of finite Galois extensions of $\Q$ with Galois groups
$G_n$ such that $G_\infty:=\prolim G_n$ is a $p$-adic Lie group of dimension at least $1$. Moreover, we assume 
that only finitely many primes ramify in $K_\infty := \bigcup_n K_n$. 

The classical example is the cyclotomic tower $K_n:=\Q(\zeta_{p^n})$ with $\zeta_{p^n}$ a $p^n$-th root of unity. A non-abelian example is
the tower  $K_n:=\Q(E[p^n])$, where $E[p^n]$ are the $p^n$-torsion points of an elliptic curve $E$
without CM defined over $\Q$.

The {\em Iwasawa algebra} 
\[ \Lambda := \Z_p[[G_\infty]]=\prolim \Z_p[G_n] \]
is the ring of $\Z_p$-valued distributions on $G_\infty$. It
is a possibly non-commutative Noetherian semi-local ring. If $G_\infty$ is in addition a pro-$p$-group without
$p$-torsion it is even a regular and local ring.

For the cyclotomic tower, $ \Lambda\isom \Z_p[G_1][ [t] ]$ is the classical Iwasawa algebra.
For the tower of $p^n$-torsion points of $E$, the Iwasawa algebra was studied 
by
Coates and Howson \cite{CH1}, \cite{CH2}.
Modules over such algebras are studied recently by Venjakob \cite{V} and by
Coates-Schneider-Sujatha \cite{CSS}.

We are concerned with the complex of  $\Lambda$-modules
$R\Gamma(\Z[1/S],\Lambda\tensor_{\Z_p}T_p(k))$ and $(\Lambda\tensor_{\Z} T_B(k-1))^+$. They are perfect complexes. Note that
\[ R\Gamma(\Z[1/S],\Lambda\tensor_{\Z_p}T_p(k))=
  \prolim R\Gamma(\Oh_{K_n}[1/S],T_p(k)) \]
where $\Oh_{K_n}$ is the ring of integers of $K_n$.

\subsection{Formulation of the non-abelian Main Conjecture}
The Main Conjecture can be viewed as a Bloch-Kato conjecture for the deformed ``motive''
$\Lambda\tensor M$ with coefficients in $\Lambda$.

Recall from \ref{conj} that the generators $\delta_p(G_n,M,k)$ are compatible under the 
transition maps $\Q_p[G_n]\to \Q_p[G_{n-1}]$. They define
\begin{multline*}
\delta_p(G_\infty,M,k)=\prolim \delta_p(G_n,M,k) \in\\
\prolim \Big[\det_{\Q_p[G_n]}R\Gamma(\Z[1/S],\Q_p[G_n]\tensor M_p(k))
\tensor\det_{\Q_p[G_n]}(\Q_p[G_n]\tensor M_B(k-1))^+
\Big]\ ,
\end{multline*}
 more precisely an element of $\prolim \Hom_{V(\Q_p[G_n])}(\1_{\Q_p[G_n]},\ \cdot\ )$.

The map $\Lambda\to \Q_p[G_n]$ induces an isomorphism
\[\Q_p[G_n]\tensor_{\Lambda}  R\Gamma\big(\Z[1/S],\Lambda\tensor T_p(k)\big)\to R\Gamma\big(\Z[1/S],\Q_p[G_n]\tensor_{\Q_p} M_p(k)\big)\ .\]

\begin{conj}[Non-abelian Main Conjecture]\label{mainconj}
Let $M$ and $S$ be as in \ref{notation}, $G_\infty$ as in \ref{3-1}, $p\neq 2$, $T_B\subset M_B$ a lattice such that $T_p:=T_B\otimes\Z_p$
is Galois stable and  $k$ big enough (cf. section \ref{2.4}). Then
$\delta_p(G_\infty,M,k)$ is induced by a generator
\[ \deltatilde_p(G_\infty,M,k)\in 
\left[\det_{\Lambda}R\Gamma(\Z[1/S],\Lambda\tensor T_p(k))\tensor
\det_{\Lambda}(\Lambda\tensor T_B(k-1))^+
\right]\ .
\]
\end{conj}
The conjecture translates into the
Iwasawa Main Conjecture in the case of Dirichlet characters or CM-elliptic curves.
See section \ref{sectionDir} for more details.\\

\rem 
a) The conjecture is independent of the choice of lattice $T_B$. The
correction factor $(\Lambda\tensor T_B(k-1))^+$ compensates
different choices of lattice.\\
b) Perrin-Riou \cite{PR} has defined a $p$-adic $L$-function and stated a Main
Conjecture for motives in the abelian case. She starts at the other 
side of the functional equation, where the exponential map of Bloch-Kato
comes into play. Her main tool is the ``logarithme {\'e}largi'', which maps
Galois cohomology over $K_\infty$ to a module of $p$-adic analytic nature. It would 
be interesting to compare her approach with the above. \\
c) A Main Conjecture for motives and the cyclotomic tower was formulated by Greenberg
\cite{Gr1}, \cite{Gr2}. Ritter and Weiss consider the case of the cyclotomic tower over a finite non-abelian
extension \cite{RitterWeiss}.

\begin{prop}[see section \ref{proof}]\label{prop}
The equivariant Bloch-Kato conjecture for $M$, $k$ and all $G_n$ 
is equivalent to the Main Conjecture for $M$, $k$ and $G_\infty$.
\end{prop}

\subsection{Twisting}
Assume that $T_p$ becomes trivial over $K_\infty$, for example let 
$G_\infty $ be  the image of $\Gal(\bar{\Q}/\Q)$ in $\Aut(T_p)$.
Let $T_p^\triv$ be the $\Z_p$-module underlying $T_p$ with trivial operation of
the Galois group. The map $g\tensor t\mapsto g\tensor g^{-1}t$
induces an equivariant isomorphism 
$\Lambda\tensor_{\Z_p} T_p\isom \Lambda\tensor_{\Z_p} T_p^\triv$.
Hence there is an isomorphism
\begin{multline*}
\det_{\Lambda}R\Gamma(\Z[1/S],\Lambda\tensor T_p(k))\tensor
\det_{\Lambda}(\Lambda\tensor T_B(k-1))^+
\isom\\
\det_{\Lambda}R\Gamma(\Z[1/S],\Lambda\tensor T_p^\triv(k))\tensor
\det_{\Lambda}(\Lambda\tensor T_B^\triv(k-1))^+
\end{multline*}
Note that $T_B^\triv$ can be viewed as a lattice in the Betti-realization
of the trivial motive $h^0(\Q)\tensor M^\triv=\Q(0)\tensor M^\triv$ where $M^\triv$ is $M_B$ considered
as $\Q$-vector space.

\begin{cor} \label{twist} If the Main Conjecture is true for $M$ and 
$\Q(0)\tensor M^\triv$ and $k$, then
\[ \deltatilde_p(G_\infty,M,k)=\deltatilde_p(G_\infty,M^\triv,k) \]
up to an element in $K_1(\Lambda)$
under the above isomorphism.
\end{cor}
\rem Even if $T_p$ is not trivial over $K_\infty$,
the same method allows to twist with a motive whose
$p$-adic realization is trivial over $K_\infty$. A particular
interesting case is the motive $\Q(1)$ if $K_\infty$ contains
the cyclotomic tower. It allows to pass from values of the $L$-function
at $k$ to values at $k+1$.\\

\noindent{\bf Strategy:}\label{strategy}
This observation allows the following strategy for proving
the Main Conjecture and the Bloch-Kato conjecture for all motives:
\begin{itemize}
\item first prove the equivariant Bloch-Kato conjecture for the motive $h^0(\Q)=\Q(0)$, 
one fixed $k$ and all finite groups $G_n$. For $k=1$ this is an equivariant
class number formula.
\item by proposition \ref{prop} this implies the Main Conjecture for the
motives $\Q(k)\tensor M^\triv$ and all $p$-adic Lie groups $G_\infty$.
\item for any
motive $M$ there is a $K_\infty$ such that $T_p$ becomes trivial. Using
corollary \ref{twist} it remains to show that $\deltatilde_p(G_\infty,M^\triv,k)$
induces $\delta_p(G_n,M,k)$ for all $n$. This is a compatibility
conjecture for elements in motivic cohomology and allows to reduce to
the case of number fields.
\item the equivariant Bloch-Kato conjecture follows by \ref{prop}. 
\end{itemize}

\section{Relation to classical Iwasawa theory in the critical case}

\subsection{Characteristic ideals}
We restrict to the case $G_\infty $ a pro-$p$-group without $p$-torsion. In this case
the Iwasawa algebra is local and Auslander regular (\cite{V}). Its total ring of quotients
is a skew field $D$. Then
$K_0(\Lambda)\isom K_0(D)\isom \Z$, $K_1(\Lambda)=(\Lambda^*)^\ab$, and
$K_1(D)=(D^*)^\ab$ where $\cdot^\ab$ denotes the abelianization of the
multiplicative group.

Let $\Th$ be the category of finitely generated $\Lambda$-torsion modules.
The localization sequence for $K$-groups implies an exact sequence
\[ (\Lambda^*)^\ab\to (D^*)^\ab\to K_0(\Th)\to 0\]
If $X$ is a $\Lambda$-torsion module, then we call its class in $K_0(\Th)$
the {\em characteristic ideal}. By the above sequence it is
an element of $D^*$ up to $[D^*,D^*]\bild \Lambda^*$. If $G_\infty$ is
abelian, $K_0(\Th)$ is nothing but the group of fractional ideals
that appears in classical Iwasawa theory.

The characteristic ideal can also be computed from the theory of
determinants. The class of $X$ in $K_0(\Lambda)$ is necessarily $0$,
hence there exists a generator $x\in \det_\Lambda(X)$. Its image
in $D\tensor\det_\Lambda(X)=\det_D(0)=\1_D$ is an element of $K_1(D)$.
This construction yields a well-defined element of 
$K_1(D)/\bild K_1(\Lambda)\isom K_0(\Th)$, 
 in fact
the inverse of the characteristic ideal of $X$.

Note that a complex is perfect if and only if it is a bounded complex with
finitely generated cohomology. 
Such complexes also have characteristic ideals if their cohomology is $\Lambda$-torsion.\\

\rem Coates, Schneider and Sujatha study the category of $\Lambda$-torsion modules
in \cite{CSS}. In particular, they also define a notion of
characteristic ideal as object of $K_0(\Th^b/\Th^1)$ where
$\Th^b/\Th^1$ denotes the quotient category of bounded finitely generated $\Lambda$-torsion modules by the sub-category of
pseudo-null modules. They construct a map
\[ K_0(\Th)\to K_0(\Th/\Th^1)\to K_0(\Th^b/\Th^1) \]
which maps the class of a module to the characteristic ideal in their sense.
If $G_\infty$ is abelian, then
the two maps are isomorphisms and all notions of characteristic ideals
agree. In the general case, we do not know whether the map is injective.
However, it seems to us that the problem is not so much in passing to
the quotient category modulo pseudo-null modules but rather in projecting
to the bounded part.

\subsection{Zeta distributions}\label{zetadistr}
Let $M$, $k$, $S$ and $G_\infty$ as before. 
Assume 
\[ \hm^1(\Z,\Q[G_n]\tensor M(k))=0 \text{ for all $G_n$.} \]
For $k$ big enough,
this implies that $M_B(k-1)^+=0$ and $K_n$ totally real.
The motives $ \Q[G_n]\tensor M(k)$ are 
 {\em critical} in the sense of Deligne. Note that the only motives expected 
to be critical and
to satisfy our condition $k$ big enough (see \ref{2.4}) are Artin motives (with $k>1)$.

In this case, the Beilinson conjecture asserts that
$L_S(G_n,M^\lor ,1-k)\in Z(\Q[G_n])^*$  (no leading coefficients has to
be taken). 
We call
\[\Lh_S(G_\infty,M^\lor,1-k)=\prolim L_S(G_n,M^\lor,1-k)    \in \prolim Z(\Q_p[G_n])^*\]
the {\em zeta distribution}.

Let $f,g\in \Lambda$ such that the images $f_n,g_n\in \Z_p[G_n]$ are units in $\Q_p[G_n]$. Via the reduced norm, they define a distribution
\[ (\rn (f_ng^{-1}_n))_n\in \prolim Z(\Q_p[G_n])^*. \]

\rem It is not clear to us if the
class of $f/g\in K_1(D)=(D^*)^\ab$ is uniquely determined by the sequence $f_ng^{-1}_n$.
In the abelian case this is true and $f/g$ is a generalization of Serre's  pseudo measure (cf. \cite{Serre}).\\

In this case the complexes $R\Gamma(\Oh_{K_n}[1/S], T_p(k))$ are torsion. Hence
the complex $R\Gamma(\Z[1/S], \Lambda\otimes T_p(k))=\prolim_n R\Gamma(\Oh_{K_n}[1/S], T_p(k))$
is bounded and its cohomology is $\Lambda$-torsion (see \cite{Ha}). The main conjecture \ref{mainconj} takes
the following form:

\begin{conj}\label{criticalconj}
Let $M$ be an Artin motive, $k>1$, $S$, $G_\infty$ as before (in particular $G_\infty$ pro-$p$ and without
$p$-torsion)
and $\Q[G_n]\tensor M(k)$ 
critical for all $n$.
There exist $f,g\in \Lambda$ such that the induced distribution $(\rn (f_ng^{-1}_n))_n\in \prolim Z(\Q_p[G_n])^*$
is the zeta distribution $\Lh_S(G_\infty,M^\lor,1-k)$ and the characteristic ideal 
\[
[R\Gamma(\Z[1/S], \Lambda\otimes T_p(k))[1]]\in K_0(\Th)
\]
coincides with the image of $fg^{-1}\in (D^*)^\ab$. 
\end{conj}

\rem a) The conjecture is isogeny invariant, i.e., independent of
the choice of lattice $T_p$. The correction term
$(\Lambda\tensor T_B(k))^+$ vanishes.\\
b) In the abelian case this means that the zeta distribution is a pseudo measure
and generates the characteristic ideal. \\
c) In the case of the cyclotomic tower, a similar conjecture is formulated by
Greenberg, \cite{Gr1}, \cite{Gr2}.\\
d) If $G_\infty$ is abelian, the above conjecture is easily seen to be implied by conjecture \ref{mainconj}.
The argument also works in the non-abelian case if the set of all elements of $\Lambda$ which, for all $n$, are
units in $\Q_p[G_n]$ is an Ore set.

\section{Examples}\label{sectionDir}
\subsection{Dirichlet characters}
Let $\chi$ be a Dirichlet character, $V(\chi)$ its associated motive
with coefficients in $E$. 
Let $\Q_\infty =\bigcup_n\Q_n$ be the cyclotomic $\Z_p$-extension of $\Q$
and $G_\infty =\Gal(\Q_\infty/ \Q)=\prolim_n G_n$. In this case the equivariant
$L$-function is $L_S(G_n,V(\chi), s)= (L_S(\rho\chi,s ))_\rho$, where $\rho$ runs through all characters of $G_n$
and $L_S(\rho\chi ,s )$ is the Dirichlet  $L$-function associated to $\rho \chi $.
Let $k$ be big enough, i.e., $k>1$.\\\\
{\bf Critical case:} $\chi(-1)=(-1)^k$.\\
Here $\hm^1(\Z, E[G_n]\otimes V(\chi)(k))=0$ for all $n$. As in section \ref{zetadistr}, the
equivariant $L$-values give rise to the zeta distribution $\Lh_S(G_\infty, V(\chi)^\lor, 1-k)\in \prolim_n E[G_n]$.
It is a classical calculation (Stickelberger elements) that this is in fact a pseudo measure,
which gives rise to the Kubota-Leopoldt $p$-adic $L$-function. Let $\Oh\subset E$ be the ring of
integers, $\Lambda=\Oh_p[ [ G_\infty] ]$ the Iwasawa algebra and $T_p(\chi)\subset V_p(\chi)$ a Galois stable lattice.
The Iwasawa Main Conjecture \ref{criticalconj} amounts to the following theorem:
\begin{thm}\label{critical} The zeta distribution $\Lh_S(G_\infty, V(\chi)^\lor, 1-k)$ generates
\[
\det^{-1}_{\Lambda} H^1(\Z[1/S], \Lambda\otimes T_p(\chi)(k))\otimes \det_\Lambda H^2(\Z[1/S], \Lambda\otimes T_p(\chi)(k)).
\]
\end{thm}
\rem This is a reformulation of the main theorem of Mazur and Wiles in \cite{Mazur-Wiles}.
There is an extension to the case of totally real fields by Wiles \cite{Wiles} and
an equivariant version by Burns and Greither  \cite{BuGr2}.
\\\\
{\bf Non-critical case:} $\chi(-1)=(-1)^{k-1}$.\\
Here $\hm^1(\Z, E[G_n]\otimes V(\chi)(k))$ has $E[G_n]$-rank $1$. It is a theorem
of Borel (resp. Soul{\'e}) that $r_\Dh\otimes \R$ (resp. $r_p\otimes \Q_p$) is an isomorphism.
By a theorem of Beilinson-Deligne (see \cite{HuWi} or \cite{HuKi-l-adic}), the image of 
 $\delta_p(G_n,V(\chi),k)$  under
$r_p$ is given by 
\[
c_k(G_n,t_p(\chi))^{-1}\otimes t_p(\chi)(k-1),
\]
where 
\[ c_k(G_n,t_p(\chi))\in H^1(\Z[1/S], \Oh_p[G_n]\otimes T_p(\chi)(k))\]
 is a twist of a cyclotomic unit and
$ t_p(\chi)(k-1)$ is a generator of $T_p(\chi)(k-1)$. Let $c_k(G_\infty ,t_p(\chi)):= \prolim_n c_k(G_n,t_p(\chi))$. 
\begin{thm}\label{non-critical} There is a canonical isomorphism of $\Lambda$-determinants
\begin{multline*}
\det_{\Lambda }\left( H^1(\Z[1/S],\Lambda \otimes T_p(\chi)(k))/c_k(G_\infty ,t_p(\chi))\right) 
\isom\\ \det_{\Lambda }H^2(\Z[1/S],\Lambda \otimes T_p(\chi)(k)).
\end{multline*}
\end{thm}
\rem  For $ p\nmid \ord(\chi)$ this is a consequence of theorem \ref{critical} and was
shown directly by Rubin \cite{Rubin} with Euler system methods. The restriction 
at the order of $\chi$ is removed in Burns-Greither \cite{BuGr} and Huber-Kings \cite{HuKi} by different methods.\\

The Tamagawa number conjecture for $V(\chi)(r)$ (and hence for $h^0(F)(r)$ with $F$ an abelian number field)
can be deduced from theorems \ref{critical} and \ref{non-critical}, see Burns-Greither \cite{BuGr} or Huber-Kings
\cite{HuKi}. Previous partial results were proved
in  Mazur-Wiles \cite{Mazur-Wiles}, Wiles \cite{Wiles}, Kato \cite{Kato-mc}, \cite{KatoGalaxy}, Kolster-Nguyen Quang Do-Fleckinger \cite{Kolster-etc} and Benois-Nguyen Quang Do\cite{BenNg}.

We would like to stress that the strategy \ref{strategy} is used in Huber-Kings \cite{HuKi}
to prove theorems \ref{critical}, \ref{non-critical} and the Tamagawa number conjecture
from the class number formula.

\subsection{Elliptic curves} 
Let $E$ be an elliptic curve over an imaginary quadratic field $K$ with CM
by $\Oh_K$.
The motive $h^1(E)$ considered with coefficients in $K$ decomposes into $V(\psi)\oplus V(\bar{\psi})$,
where $\psi$ is the Gr{\"o}ssencharacter associated to $E$. The $L$-function of $V(\psi)$
is the Hecke $L$-function of $\psi$, which has a zero of order $1$ at $2-k$, where
$k\geq 2$. Let $S=Np$, where $N$ is the conductor of $\psi$ and let $K_n:=K(E[p^n])$.
 
It is not known if $\hm^1(\Oh_K,K[G_n]\otimes V(\psi)(k))$ has $K[G_n]$-rank $1$
but Deninger \cite{De}  shows that $r_\Dh \otimes \R$ is surjective
and that the Beilinson conjecture holds.
It is a result of Kings \cite{Ki} that the image in {\'e}tale cohomology of the zeta element
 $\delta_p(G_n,V(\psi),2-k)$ given by Beilinson's Eisenstein symbol is
given by 
\[e_k(G_n, t_p(\psi))^{-1}\otimes t_p(\psi)\ ,\] 
where
 $e_k(G_n, t_p(\psi))\in  H^1(\Z[1/S],\Oh_p[G_n]\otimes T_p(\psi )(k))$ is the twist of an  elliptic unit.
Let $\Lambda:=\Oh_p[ [G_\infty] ]$ and $e_k(G_\infty ,t_p(\psi ))=\prolim_n e_k(G_n, t_p(\psi))$.
\begin{thm} 
There is a canonical isomorphism of determinants
\begin{multline*}
\det_{\Lambda }\left( H^1(\Z[1/S],\Lambda \otimes T_p(\psi )(k))/e_k(G_\infty ,t_p(\psi ))\right) 
\isom\\ \det_{\Lambda }H^2(\Z[1/S],\Lambda \otimes T_p(\psi )(k)).
\end{multline*}
\end{thm}

\rem 1) This is a reformulation of Rubin's Iwasawa Main Conjecture \cite{Ru}.\\
2) In \cite{Ki} the (absolute) Bloch-Kato
conjecture for $V(\psi)$ is deduced from this under the condition that $H^2(\Z[1/S],T_p(\psi )(k))$ is finite
(fulfilled for almost all $k$ for fixed $p$).\\

Kato \cite{Ka} has investigated the case of elliptic curves over $\Q$ and the cyclotomic tower.
His approach to the Birch-Swinnerton-Dyer conjecture uses the idea of twisting 
cup-products of Eisenstein symbols to the value of the $L$-function at $1$.
As a consequence he can prove one inclusion of the Iwasawa main conjecture in this case.
 The result supports our general philosophy of
twisting to the case of number fields. \\

\section{Proof of proposition \ref{prop}}\label{proof}
We want to give the proof of proposition \ref{prop}. The implication
from the Main Conjecture to the equivariant Bloch-Kato conjecture
is trivial. Conversely, we have to show the following abstract statement:

\begin{lemma}
Let $\nabla\in V(\Lambda)$ and 
$\deltatilde(n)\in\Z_p[G_n]\tensor \nabla$ generators such that their
images $\delta(n)\in \Q_p[G_n]\tensor \nabla$ are compatible
under transition maps. Then there is a generator $\deltatilde'(\infty)\in \nabla$ inducing
all $\delta(n)$.
\end{lemma}

The proposition follows with  $\deltatilde(n)=\deltatilde_p(G_n,M,k)$ and
\[ \nabla=\det_{\Lambda}R\Gamma(\Z[1/pS],\Lambda\tensor T_p(k))\tensor
\det_{\Lambda}(\Lambda\tensor T_B(k-1))^+
\]

We now prove the lemma.
We first reduce to a statement about elements of $K_1$. 
By assumption,
$\Z_p[G_n]\tensor\nabla $ has a generator, in particular, its isomorphism class is
zero in $K_0(\Z_p[G_n])$.
As $K_0(\Lambda)\to \prolim K_0(\Z_p[G_n])$ is an isomorphism, this
implies that the class of $\nabla$ is zero in
$K_0(\Lambda)$. Without loss of generality we can assume $\nabla=\1_\Lambda$.
Recall that by our convention, a generator of $\1_A$ is nothing but an element of
the abelian group $K_1(A)$ for all rings $A$.

Let 
$B_n=\bild K_1(\Z_p[G_n])\to K_1(\Q_p[G_n])$. By assumption $\delta(n)\in B_n$.
There is a system of short exact sequences
\[ 0\to SK_1(\Z_p[G_n])\to K_1(\Z_p[G_n])\to B_n\to 0 \]
By \cite{CurtisReiner} 45.22 the groups $SK_1(\Z_p[G_n])$ are finite.
The system of these groups is automatically Mittag-Leffler. Hence we get a
surjective map
\[ \prolim K_1(\Z_p[G_n])\to \prolim B_n \]
The system $(\delta(n))_n$ has a preimage 
$(\deltatilde'(n))_n\in \prolim K_1(\Z_p[G_n])$.

All $\Z_p[G_n]$ are semi-local, hence by \cite{CurtisReiner} 40.44
\[ K_1(\Z_p[G_n])\isom \Gl_2(\Z_p[G_n])/E_2(\Z_p[G_n]) \]
where $E_2$ is the subgroup of elementary matrices. We represent
$\deltatilde'_p(n)$ by an element of $\Gl_2(\Z_p[G_n])$. By assumption
the image of $\deltatilde'(n)$ 
in $K_1(\Z_p[G_{n-1}])$
differs from $\deltatilde'(n-1)$ by some
elementary matrix in $E_2(\Z_p[G_{n-1}])$. Elementary matrices can be lifted to elementary matrices
in $\Gl_2(\Z_p[G_n])$. Hence we can assume that the elements
$\deltatilde'(n))\in \Gl_2(\Z_p[G_n])$ form a projective system.
The system defines an element 
\[ \deltatilde'_p(n)\in \Gl_2(\Lambda)\]
whose class in $K_1(\Lambda)$ has the necessary properties.

\end{document}